\documentclass[a4paper,10pt]{article}
    	\usepackage{enumerate}
    	\usepackage{color}
    	\usepackage[applemac]{inputenc}
    	\usepackage[english]{babel}
    	\usepackage[T1]{fontenc}
    	\usepackage{graphicx}
    	\usepackage{amsfonts,amssymb,amsmath,latexsym,amsthm}
    	\usepackage{textcomp}
    	\usepackage[pdftex]{hyperref}
    	\usepackage{geometry}
    	\usepackage{tikz}
     \usepackage{bm}

    	\author{L. G. Martins$^1$,
     M. V. Flamarion$^2$,  R. Ribeiro-Jr$^3$}
     
    	\title{A forced Boussinesq model with  a sponge layer}

    	\date{}
    	
\begin{document}
    	\maketitle

    	\begin{center}

    		{\footnotesize $^1$ UFPR/Federal University of Paran‡ \\
    	luiz.martins1@ufpr.br }

    		\vspace{0.3cm}
    
    		{\footnotesize $^{2}$
      Unidade Acad{\^ e}mica do Cabo de Santo Agostinho, 
    	UFRPE/Rural Federal University of Pernambuco, BR 101 Sul, Cabo de Santo Agostinho-PE, Brazil,  54503-900 \\
    	marcelo.flamarion@ufrpe.br}
    
    		\vspace{0.3cm}

    		{\footnotesize $^3$
      UFPR/Federal University of Paran\'a,  Departamento de Matem\'atica, Centro Polit\'ecnico, Jardim das Am\'ericas, Caixa Postal 19081, Curitiba, PR, 81531-980, Brazil \\ robertoribeiro@ufpr.br}
    	
    	

    	\end{center}

    	
    	\begin{abstract} 
    	\noindent The movement of water waves is a topic of interest to researchers from different areas. While their propagation is described by Euler equations, there are instances where simplified models can also provide accurate approximations. A well-known reduced model employed to study the wave dynamics is the Boussinesq model. Despite being extensively studied, to our knowledge, there is no research available on a Boussinesq model featuring a sponge layer. Therefore, in this work, we present a Boussinesq model with a sponge layer. Furthermore, we carry out a numerical investigation to explore the advantages and limitations of the proposed model. For this purpose, we compare the numerical solutions of the model with and without the sponge in three different scenarios. The numerical solutions are computed by a pseudospectral method. Our results show that the Boussinesq model with a sponge layer is numerically stable and advantageous because it is able to absorb low-amplitude waves, allowing it to run the numerical simulations for long periods of time without requiring a large spatial domain, but it is not able to absorb high-amplitude waves. 
    		\end{abstract}

    		\section{Introduction}
    
    The study of water waves is a subject that catches the attention of researchers from different areas, such as engineers, physicists, and mathematicians. Describing their behavior may be a difficult task because their motion can be influenced by several factors, such as the shape of the seabed, the wind and the presence of a current. A better understanding of the motion of the water waves can provide, for instance, more security and confiability in designing coastal and offshore structures  \cite{MejiaGrajales:2021}.
    
     Research on the propagation of water waves over topographic obstacles \cite{Baines, Pratt, Kim, Lee, LeeWhang, Paul}, ship waves, and ocean waves generated by storms \cite{Johnson, Wu1, Wu2} has been a focal point in applied mathematics. A particular focus has been on trapped waves, which are waves that persistently remain confined within specific spatial regions, often around topographic obstacles or under the influence of external forces. Grimshaw et al. \cite{Grimshaw94} delved into the interaction of solitary waves with a small-amplitude external force using the forced Korteweg-de Vries (fKdV) equation, exploring both asymptotic and numerical approaches. Their investigations revealed regimes in which solitary waves either partially or entirely remained trapped by the external force. The work of Lee and Whang \cite{LeeWhang} and Lee \cite{Lee} extended this understanding by examining a scenario involving a two-bumped obstacle. They derived solutions for the fKdV equation, demonstrating the persistence of trapped waves between the two obstacles over a specific duration. Additionally, they conducted a thorough analysis of the numerical stability of these trapped waves by perturbing their initial amplitudes and the topography height. In a similar vein, Kim and Choi \cite{Kim} contributed to the field by establishing that trapped waves must overcome a distinct energy barrier to escape the region between the obstacles. This finding adds a crucial dimension to the understanding of the dynamics of trapped waves, shedding light on the intricate interplay between external forces, topographic features, and the persistence of wave confinement.

    	In addition to the forced Korteweg-de Vries equation, another well-known reduced model employed in the study of wave dynamics is the Boussinesq model. Although extensively studied, to our knowledge, there is no research available on a Boussinesq model featuring a sponge layer.  From a physical perspective, the sponge can be seen as a damping term placed at the boundary of the computational domain 
    	 that absorbs small-amplitude waves and let them leave the domain without reflecting any significant signature inside.  This approach has demonstrated successful application in the design of numerical formulation   in other models. For example, in computational aeroacoustic \cite{Mani:2012}, in the fKdV equation \cite{GrimShawMalewoong:2019},  in  the Ostrovysky equation  \cite{Alias:2019} and in  the forced Whitham equation \cite{Flamarion:sd}.

    Therefore, to fill this gap, we present a Boussinesq model with a sponge layer. Our model  is an adaptation of the Boussinesq model proposed by Chen \cite{Chen:2003}, originally formulated for waves propagating on a channel with a moving bottom and under the effect of a variable pressure on the free surface.
    Besides, a numerical investigation is carried out in order to explore  the advantages and limitations of  the model that we proposed.
      For this purpose, we compare the numerical solutions of the model with and without the sponge in the following scenarios: i) propagation of a traveling solitary wave in a channel with an even bottom and a constant pressure on the free surface; ii) study of trapped waves; iii) waves generated by a current-topography interaction.   The numerical solutions are computed by a pseudospectral method. Our results show that the Boussinesq model with a sponge layer is advantageous because it is able to absorb low-amplitude waves, allowing it to run the numerical simulations for long periods of time without requiring a large spatial domain, but it is not able to absorb high-amplitude waves.

    The paper is structured as follows. The mathematical formulation and numerical scheme are presented in section 2 and 3.  The results are described in section 4. The final remarks are made in section 5.

    \section{Mathematical Formulation}

    \subsection{Boussinesq Equations}
    
    Consider an inviscid, incompreesible, and irrotational fluid with a constant density $\rho$ in a two-dimensional space $(x,y)$. We introduce a Cartesian coordinate system $(x,y)$ with gravity pointing in the negative $y$-direction and $y = 0$ being the undisturbed free-surface. The fluid is under the effect of a variable pressure on the free surface $P(x,t)$ that varies in both space and time. The bottom boundary of the fluid is at $y = -1 + h(x,t)$ and the upper boundary is free to move and denoted by $\eta(x,t)$. As derived by Chen \cite{Chen:2003}, the Boussinesq equations in dimensionless variables read as
    \begin{equation}
    \label{equacoesdeBoussinesqpadrao}
        \begin{cases}
        \eta_t + ((1-\alpha h + \alpha \eta)u)_x = h_t,\\
        u_t + \eta_x + \alpha uu_x - \dfrac{\beta}{3}u_{xxt} = -\dfrac{\beta}{2}h_{xtt} - \alpha P_x.
        \end{cases}
    \end{equation}
    In the absence of external forces ($P=0=h$), the system (\ref{equacoesdeBoussinesqpadrao}) supports left-going and right-going solitary waves that closely resemble KdV soliton solutions up to $\mathcal{O}(\alpha^{2},\alpha\beta)$, as elaborated in \cite{Whitham:1974}. Moving forward, we introduce the Boussinesq equations with a sponge layer.

    \subsection{Inserting a sponge layer}
    System \eqref{equacoesdeBoussinesqpadrao} can been implemented numerically with periodic boundary conditions applying the Fourier Transform. It is expected to obtain an accurate approximation as long as the solution does not approach the boundaries \cite{Shentangwang:2011}. When it reaches the boundaries, due to the periodic conditions, the waves re-enter from them and start interacting with the solution.

\begin{figure}[!ht]
            	\centering
            	\includegraphics{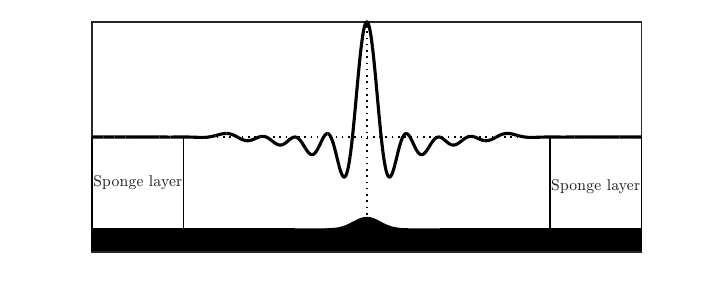}
        \caption{Sponge scheme}
        \label{esquema-esponja}
            \end{figure}

    Although extending the spatial domain seems a practical solution for simulating the numerical experiments for long periods of time, it increases the computational cost. A possible alternative is inserting an artificial sponge near the domain boundaries, which is a common practice due to its simplicity, efficiency and non-stiff nature \cite{Mani:2012}. Physically, the sponge acts as an damping system that absorbs waves when they get close to the boundaries, as showed in the schematic figure \ref{esquema-esponja}.
    
    Mani \cite{Mani:2012} has provided practical guidelines for the development of sponge layers. To show that they can achieve the desired accuracy, the author presented a convective vortex/sponge interaction using the Euler model. Additionally, other examples of inserting a sponge can be found in Grimshaw and Malewoong \cite{GrimShawMalewoong:2019}, who inserted a sponge in the fKdV equation, Alias et al. \cite{Alias:2019}, who combined a sponge with other techniques to numerically study the Ostrovysky equation, and Flamarion et al. \cite{Flamarion:sd}, who coupled a sponge in the forced Whitham equation.
    
    We are interested in inserting a sponge in the Boussinesq equations, so we are going to consider a particular case of the system \eqref{equacoesdeBoussinesqpadrao}. In what follows we present the motivation for the sponge layer considered in this paper. By considering $h = 0, P = 0$, and $\alpha,\beta \rightarrow 0^+$, we obtain
    \begin{equation}
    \label{EQBOUSSREDZ}
    \begin{cases}
    \eta_t + u_x = 0, \\  
    u_t + \eta_x  = 0.
    \end{cases}
    \end{equation} 
    
    The system \eqref{EQBOUSSREDZ} can be simplified to a single equation by differentiating the first equation with respect to $x$, differentiating the second equation with respect to $t$, and then taking the difference. Hence, we get
    \begin{equation*}
     \eta_{tt} - \eta_{xx} = 0,
    \end{equation*}
    which is the wave equation in the adimensional variables. Adding the initial conditions, we arrive at the following problem:
    \begin{equation}
    \label{edpondas}
    \begin{cases}
    \eta_{tt} - \eta_{xx} = 0,\\
    \eta(x,0) = f(x), \\
    \eta_t(x,0) = g(x),
    \end{cases}
    \end{equation}
    whose solution is
    \begin{equation}
    \label{eqondalinearsemdispersao}
        \eta(x,t) = \dfrac{f(x+t)+f(x-t)}{2} + \dfrac{1}{2}\int_{x-t}^{x+t} g(s) \ ds.
    \end{equation}
    
    A damping system can be inserted in problem \eqref{eqondalinearsemdispersao} by modifying the PDE. Thus, we have
    \begin{equation}
    \label{edpondaamortecimento0}
    \begin{cases}
    \eta_{tt} - \eta_{xx} + 2b\eta_t + b^2\eta = 0, \\
    \eta(x,0) = f(x), \\
    \eta_t(x,0) = g(x),
    \end{cases}
    \end{equation}
    where $b > 0$.
    
\begin{figure}[!htb]
        \centering
        \includegraphics{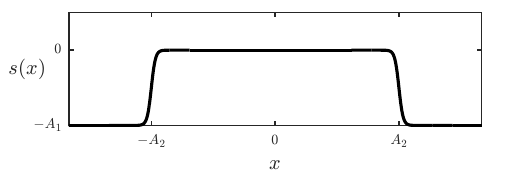}
        \caption{Sponge sketch}
        \label{spongegraph}
    \end{figure}

    The solution of the new problem \eqref{edpondaamortecimento0} is given by
    \begin{equation}
        \eta(x,t) = e^{-bt}\left(\dfrac{f(x+t)+f(x-t)}{2} + \dfrac{1}{2}\int_{x-t}^{x+t} g(s) \ ds + \dfrac{b}{2}\int_{x-t}^{x+t} f(s) \ ds  \right).
        \label{newsolution}
    \end{equation}
    Notice that for small values of $b$, the solution \eqref{newsolution} is nearly identical to the solution \eqref{eqondalinearsemdispersao}. In addition, for large values of $t$, we have $\eta(x,t)\rightarrow 0$. Therefore, we would like to insert a sponge that is a non-zero constant near the boundaries and almost zero elsewhere. A function $s(x)$ that satisfies those features is
    \begin{equation}
        s(x) = \dfrac{A_1}{2}\biggl(\tanh(x+A_2) - \tanh(x-A_2)\biggr) - A_1, \label{esponjaa}
    \end{equation}
    with $A_1,A_2 > 0$. The parameter $A_1$ controls the intensity of the absorption, while parameter $A_2$ determines the region where the sponge has no effect. A sketch of the sponge function $s(x)$ is illustrated in figure \ref{spongegraph}.
    
    Therefore setting $b = -s(x)$, we obtain 
    \begin{equation}
        \eta(x,t) = e^{s(x)t}\left(\dfrac{f(x+t)+f(x-t)}{2} + \dfrac{1}{2}\int_{x-t}^{x+t} g(s) \ ds - \dfrac{s(x)}{2}\int_{x-t}^{x+t} f(s) \ ds  \right).
        \label{soldeondascomamortecimento}
    \end{equation}
    
    Consider the modified version of the equation  \eqref{EQBOUSSREDZ} 
    \begin{equation}
    \label{EQBOUSSREDZmod}
        \begin{cases}
        \eta_t + u_x + G_1 = 0, \\
        u_t + \eta_x - G_2 = 0,
        \end{cases}
    \end{equation}
    where $G_1 = G_1(x,t)$ and $G_2 = G_2(x,t)$ are functions to be determined. The idea is to rewrite the system \eqref{EQBOUSSREDZmod} as a single equation, and compare it with equation \eqref{edpondaamortecimento0} to determine $G_1$ and $G_2$.
    
    Rewriting the system \eqref{EQBOUSSREDZmod}, we have
    \begin{equation*}
        \eta_{tt} - \eta_{xx} + G_{1,t} + G_{2,x} = 0.
    \end{equation*}
    
    A possible choice is
    \begin{equation*}
        G_{1,t}(x,t) = -2s(x)\eta_t(x,t) \qquad \text{and} \qquad G_{2,x}(x,t) = s^2(x)\eta(x,t)
    \end{equation*}
    which implies
    \begin{equation*}
    G_1(x,t) = -2s(x)\eta(x,t) \quad \text{and} \quad G_2(x,t) = \int s^2(x)\eta(x,t) \ dx,    
    \end{equation*}
    where $G_2(x,t)$ should be understand as an antiderivative. It is important to mention that we will refer to this antiderivative as $G$ instead of $G_2$. 
    
    Therefore, the Boussinesq equations with a sponge layer are given by
    \begin{equation}
    \begin{cases}
    \eta_t + ((1-\alpha h + \alpha\eta)u)_{x} - 2s(x)\eta = h_{t}, \\
    u_t + \eta_x + \alpha uu_{x} - \dfrac{\beta}{3}u_{xxt} - G(x,t) = -\dfrac{\beta}{2}h_{xxt} - \alpha P_x,
    \end{cases}
    \label{boussinesqcomesponjaa}
    \end{equation}
    where  $G(x,t) = \int s^2(x)\eta(x,t) \ dx$.

    \section{Numerical scheme}\label{ns}
    
    \subsection{The time-dependent Boussinesq equations}
    
    Suppose that the topography and the variable pressure on the free surface move horizontally at a constant speed, that is, $h(x,t) = h(x+Ft)$ and $P(x,t) = P(x+Ft)$. By rewriting the system \eqref{boussinesqcomesponjaa} in the moving frame  $x' = x + Ft$  and abandon the notation \; $'$ \; we obtain the following system
    \begin{equation}
    \begin{cases}\label{boussinesqcomesponjaa2}
    \eta_t + F\eta_x + ((1-\alpha h + \alpha\eta)u)_{x} - 2s(x)\eta = Fh_{x}, \\  
    u_t + Fu_x + \eta_x + \alpha uu_{x} - \dfrac{\beta}{3}u_{xxt} - \dfrac{\beta}{3}Fu_{xxx} - G = -\dfrac{\beta}{2}F^2h_{xxx} - \alpha P_x.
    \end{cases}
    \end{equation}
    This system  can be interpreted as a Boussinesq model to investigate current-topography interactions with the effect of a variable pressure on the free surface. Here, $F$ denotes the Froude number.

    The numerical scheme  to solve system \eqref{boussinesqcomesponjaa2} is based on a pseudospectral method proposed by Trefethen \cite{Trefethen:2001}. We consider the system \eqref{boussinesqcomesponjaa2} with initial conditions $\eta(x,0) = \eta_o(x)$ and $u(x,0) = u_o(x)$, and we impose that the solutions, the topography, and the variable pressure on the surface satisfy the following conditions
    \begin{equation*}
        u(x,0) \rightarrow 0, \; \eta(x,0) \rightarrow 0, \; h(x,t) \rightarrow 0  \; \text{and} \; P(x,t) \rightarrow 0,
    \end{equation*}
    as $|x| \rightarrow \infty$. Thus, instead of solving the problem over $\mathbb{R}$, we can restrict on the interval $(-L,L)$, $L > 0$ sufficiently large, and we can approximate the boundary conditions by periodic boundary conditions, which will allow us to use the Fourier Transform.

    Applying the Fourier Transform in system \eqref{boussinesqcomesponjaa2} and making a few algebraic manipulations, one obtains
    \begin{equation}
    \label{eqboussinesqnova4}
        \begin{cases}
       \hat{\eta}_t = -(ik)F\hat{\eta} - (ik)\hat{u} + (ik)\alpha\widehat{hu} - (ik)\alpha\widehat{\eta u} + 2\widehat{s(x)\eta} + (ik)F\hat{h}, \\
        \hat{u}_t = \dfrac{-F(ik)\hat{u} -(ik)\hat{\eta} - \dfrac{\alpha}{2}(ik)\widehat{u^2} + \dfrac{\beta}{3}F(ik)^3\hat{u} + \dfrac{\widehat{s^2(x)\eta}}{ik} - \dfrac{\beta}{2}F ^2(ik)^3\hat{h} - \alpha (ik)\hat{P}}{1-\dfrac{\beta}{3}(ik)^2}, \\
        \hat{\eta}(k,0) = \hat{\eta}_o(k), \\
        \hat{u}(k,0) = \hat{u}_o(k).
        \end{cases}
    \end{equation}

    Now, we are going to describe the Fourier approximation of the system \eqref{eqboussinesqnova4} in MATLAB. The system is discretized on a uniform grid in the variable $x$
    \begin{equation}
        x_j = -L + j\Delta x, j = 0,1,\dots,N-1, \Delta x = \dfrac{2L}{N},
        \label{discretizacao1}
    \end{equation}
    and the functions are evaluated on the grid points, being denoted by $\eta_j = \eta(x_j,t), u_j = u(x_j,t), h_j = h(x_j,t)$, and $P_j = P(x_j,t)$. Therefore, we can define the vectors
    \begin{equation*}
        \vec{\eta} = (\eta_0,\dots,\eta_{N-1}), \quad \vec{u}=(u_0,\dots,u_{N-1}), \quad \vec{h} = (h_0,\dots,h_{N-1}) \quad \text{and} \quad \vec{P} = (P_0,\dots,P_{N-1}).
        \label{discretizacao2}
    \end{equation*}
    
    Besides, let $\bm{k}$ be the frequency vector, given by
    \begin{equation}
        \bm{k} = \left(0,1,\dots,\dfrac{N}{2}-1,0,-\dfrac{N}{2}+1,\dots,-1 \right)
        \label{discretizacao3}
    \end{equation}
    and denote by
     \begin{align}
        \tilde{\bm{\eta}} = \text{fft}(\vec{\eta}), \quad \tilde{\bm{u}} = \text{fft}(\vec{u}), \quad \tilde{\bm{h}} = \text{fft}(\vec{h}) \quad \text{and} \quad \tilde{\bm{P}} = \text{fft}(\vec{P}),
        \label{discretizacao4}
    \end{align}
    the discrete Fourier Transform of the vectors $\vec{\eta},\vec{u},\vec{h}$, and $\vec{P}$.

    Therefore, the Fourier approximation of system \eqref{boussinesqcomesponjaa2}, which is based on \eqref{eqboussinesqnova4}, is 
    \begin{equation}
    \label{eqboussinesqnova5}
        \begin{cases}
       \dfrac{d\tilde{\bm{\eta}}}{dt} = -i\bm{k}.\ast F .\tilde{\bm{\eta}} - i\bm{k}.\ast\tilde{\bm{u}} + i\bm{k}.\ast \alpha.\text{fft}(\text{ifft}(\tilde{\bm{h}}).\ast\text{ifft}(\tilde{\bm{u}})) - i\bm{k}.\alpha.\ast\text{fft}(\text{ifft}(\tilde{\bm{\eta}}).\ast\text{ifft}(\tilde{\bm{u}})) + i\bm{k}.F.\ast\tilde{\bm{h}} \\ 
       \hspace{1cm} + 2\text{fft}(\tilde{\bm{s}}.\ast\text{ifft}(\tilde{\bm{\eta}})), \\
        \dfrac{d\tilde{\bm{u}}}{dt} = \dfrac{-i\bm{k}.\ast F.\tilde{\bm{u}} -i\bm{k}.\ast\tilde{\bm{\eta}} - \dfrac{\alpha}{2}. i\bm{k}.\ast\text{fft}{((\text{ifft}(\tilde{\bm{u}}))^2)} + \dfrac{\beta}{3}.F(i\bm{k})^3.\ast\tilde{\bm{u}} - \dfrac{\beta}{2}.F^2.(i\bm{k})^3.\ast\tilde{\bm{h}} - \alpha .i\bm{k}.\ast\tilde{\bm{P}}}{1-\dfrac{\beta}{3}(ik)^2} \\
        \hspace{1cm} + \dfrac{ \dfrac{1}{ik}.\ast\text{fft}(\tilde{\bm{s}}^2.\ast\text{ifft}(\tilde{\bm{\eta}}))}{1-\dfrac{\beta}{3}(ik)^2}, \\
        \tilde{\bm{\eta}}(0) = \text{fft}(\eta_o), \\
        \tilde{\bm{u}}(0) = \text{fft}(u_o),
         \end{cases}
    \end{equation}
    with $.\ast$ denoting the component by component multiplication and \text{ifft} the inverse Fourier Transform. The equation \eqref{eqboussinesqnova5} is integrated in time using the fourth-order Runge-Kutta method.  The numerical code is available in the appendix      \ref{AppendixB}.

    \subsection{The stationary Boussinesq equations}
    
    Let us consider the system of equations \eqref{boussinesqcomesponjaa2} again. If we assume that the topography is flat and the pressure on the free surface is constant, then the system will simplify to
    \begin{equation}
    \begin{cases}
        \eta_t + F\eta_x + u_x + \alpha(\eta u)_x = 0, \\
        u_t + F u_x + \eta_x + \dfrac{\alpha}{2} (u^2)_x - \dfrac{\beta}{3} u_{xxt} - \dfrac{\beta}{3}F u_{xxx} = 0.
    \end{cases}  
    \label{sistemabase}
    \end{equation}
    
    We are interested in calculating stationary solutions, that is, we want to find functions $\eta$ and $u$ that are solution of \eqref{sistemabase} such that $\eta_t(x,t) = 0$ and $u_t(x,t) = 0, \forall t \geq 0$. For this purpose, we make $\eta_t = u_t = 0, \forall t \geq 0$ in the system \eqref{sistemabase}. Thus, we obtain that
    \begin{equation}
    \begin{cases}
        u_x + F\eta_x  + \alpha(\eta u)_x = 0, \\
        \eta_x + F u_x + \dfrac{\alpha}{2} (u^2)_x - \dfrac{\beta}{3}F u_{xxx} = 0.
    \end{cases}    
    \label{sistemabase2}
    \end{equation}
    
    For the same reason as in the previous subsection, we impose that $\eta(x,t) \rightarrow 0$ and $u(x,t) \rightarrow 0$ when $|x| \rightarrow\infty$. Hence, taking the Fourier transform on the variable $x$ in the system \eqref{sistemabase2}, we arrive at the following system
    \begin{equation}
        \begin{cases}
            \hat{u} + F\hat{\eta} + \alpha\widehat{\eta u} = 0, \\
            \hat{\eta} + F\hat{u} + \dfrac{\alpha}{2}\widehat{u^2} - \dfrac{\beta}{3}F(ik)^2\hat{u} = 0.
        \end{cases}
        \label{quaseoesquema}
    \end{equation}
    
    We will now present the Fourier approximation of the equations \eqref{sistemabase2} in MATLAB. This process is similar to the one explained in the previous subsection. The system \eqref{sistemabase2} is discretized on a uniform grid in the variable $x$, as described in the expression \eqref{discretizacao1}, we define the vectors $\vec{\eta},\vec{u},\vec{h}$ and $\vec{P}$, whose components correspond to the functions $\eta,u,h,P$ applied to the grid points. Then we define the Fourier Transforms of these vectors, according to equation \eqref{discretizacao4}. Thus, we obtain that  
    \begin{equation}
        \begin{cases}
            \tilde{\bm{u}} + F.\tilde{\bm{\eta}} + \alpha.\text{fft}(\text{ifft}(\tilde{\bm{\eta}}.\ast\tilde{\bm{u}})) = 0, \\
            \tilde{\bm{\eta}} + F\tilde{\bm{u}} + \dfrac{\alpha}{2}.\text{fft}{((\text{ifft}(\tilde{\bm{u}}))^2)} - \dfrac{\beta}{3}.F.(i\bm{k})^2.\ast\tilde{\bm{u}} = 0.
        \end{cases}
        \label{quaseoesquema2}
    \end{equation}
    
    Note that the system of equations \eqref{quaseoesquema2} has $2N+1$ unknowns ($\tilde{\eta}_0,\dots,\tilde{\eta}_{N-1},\tilde{u}_0,\dots,\tilde{u}_{N-1},F$) and $2N$ equations. For the sake of convenience, we impose that $\eta$ and $u$ are even functions, which results in
    \begin{align*}
        &\tilde{\bm{\eta}} = (\tilde{\eta}_0, \tilde{\eta}_1, \dots, \tilde{\eta}_{\frac{N}{2}}, \tilde{\eta}_{\frac{N}{2}+1},\dots,\tilde{\eta}_{N-1}) = (\tilde{\eta}_0, \tilde{\eta}_1, \dots, \tilde{\eta}_{\frac{N}{2}}, \tilde{\eta}_{\frac{N}{2}-1},\dots,\tilde{\eta}_{1}), \\
        &\tilde{\bm{u}} = (\tilde{u}_0, \tilde{u}_1, \dots, \tilde{u}_{\frac{N}{2}}, \tilde{u}_{\frac{N}{2}+1},\dots,\tilde{u}_{N-1}) = (\tilde{u}_0, \tilde{u}_1, \dots, \tilde{u}_{\frac{N}{2}}, \tilde{u}_{\frac{N}{2}-1},\dots,\tilde{u}_{1}),
    \end{align*}
    reducing the system to $N + 3$ unknowns ($\tilde{\eta}_0,\dots,\tilde{\eta}_{\frac{N}{2}},\tilde{u}_0,\dots,\tilde{u}_{\frac{N}{2}},F$) and $N + 2$ equations. 
    
    To match the number of unknowns and equations, making it possible to apply Newton's method, let's impose the following condition 
    \begin{equation*}
        \tilde{\eta}_{\frac{N}{2}} - A = 0,
    \end{equation*}
    where $A$ is the amplitude of the wave. Thus, we have the following system
    \begin{equation}
        \begin{cases}
        \tilde{\bm{u}} + F.\tilde{\bm{\eta}} + \alpha.\text{fft}(\text{ifft}(\tilde{\bm{\eta}}.\ast\tilde{\bm{u}})) = 0,  \\    
        \tilde{\bm{\eta}} + F\tilde{\bm{u}} + \dfrac{\alpha}{2}.\text{fft}{((\text{ifft}(\tilde{\bm{u}}))^2)} - \dfrac{\beta}{3}.F.(i\bm{k})^2.\ast\tilde{\bm{u}} = 0 \\
         \tilde{\eta}_{\frac{N}{2}} - A = 0.      
        \end{cases}
        \label{esquema}
    \end{equation}
    
    The system of equations \eqref{esquema} is solved using Newton's method. We call this system $G: \mathbb{R}^{n+3} \rightarrow \mathbb{R}^{n+3}$, $(\tilde{\eta}_0,\dots,\tilde{\eta}_{N/2},\tilde{u}_0,\dots,\tilde{u}_{N/2}, F) \mapsto (G_1^0,\dots,G_1^{N/2},G_2^0,\dots,G_2^{N/2},G_3)$, where the components are given by
    \begin{equation*}
        \begin{cases}
        G_1^n(\tilde{\eta}_0,\dots,\tilde{\eta}_{N/2},\tilde{u}_0,\dots,\tilde{u}_{N/2},F) = \tilde{\bm{u}} + F.\tilde{\bm{\eta}} + \alpha.\text{fft}(\text{ifft}(\tilde{\bm{\eta}}.\ast\tilde{\bm{u}})) = 0, \\
        G_2^n(\tilde{\eta}_0,\dots,\tilde{\eta}_{N/2},\tilde{u}_0,\dots,\tilde{u}_{N/2},F) = \tilde{\bm{\eta}} + F\tilde{\bm{u}} + \dfrac{\alpha}{2}.\text{fft}{((\text{ifft}(\tilde{\bm{u}}))^2)} - \dfrac{\beta}{3}.F.(i\bm{k})^2.\ast\tilde{\bm{u}} = 0, \\
        G_3(\tilde{\eta}_0,\dots,\tilde{\eta}_{N/2},\tilde{u}_0,\dots,\tilde{u}_{N/2},F) = \tilde{\eta}_{\frac{N}{2}} - A = 0,
        \end{cases}
         \label{esquema2}
    \end{equation*}
    with $n = 0,1,\dots,\dfrac{N}{2}$.
    
    The Jacobian matrix is of the form
    \begin{align*}
    \begin{bmatrix}
        \dfrac{\partial G_1^0}{\partial \tilde{\eta}_0} & \dfrac{\partial G_1^0}{\partial \tilde{\eta}_1} & \cdots & \dfrac{\partial G_1^0}{\partial \tilde{\eta}_{N/2}} & \dfrac{\partial G_1^0}{\partial \tilde{u}_0} & \dfrac{\partial G_1^0}{\partial \tilde{u}_1} & \cdots & \dfrac{\partial G_1^0}{\partial \tilde{u}_{N/2}} & \dfrac{G_1^0}{\partial F} \\
        \vdots & \vdots & \ddots &  \vdots &  \vdots &   \vdots & \ddots & \vdots & \vdots \\
        \dfrac{\partial G_1^{N/2}}{\partial \tilde{\eta}_0} & \dfrac{\partial G_1^{N/2}}{\partial \tilde{\eta}_1} & \cdots & \dfrac{\partial G_1^{N/2}}{\partial \tilde{\eta}_{N/2}} & \dfrac{\partial G_1^{N/2}}{\partial \tilde{u}_0} & \dfrac{\partial G_1^{N/2}}{\partial \tilde{u}_1} & \cdots & \dfrac{\partial G_1^{N/2}}{\partial \tilde{u}_{N/2}} & \dfrac{G_1^{N/2}}{\partial F} \\
        \dfrac{\partial G_2^0}{\partial \tilde{\eta}_0} & \dfrac{\partial G_2^0}{\partial \tilde{\eta}_1} & \cdots & \dfrac{\partial G_2^0}{\partial \tilde{\eta}_{N/2}} & \dfrac{\partial G_2^0}{\partial \tilde{u}_0} & \dfrac{\partial G_2^0}{\partial \tilde{u}_1} & \cdots & \dfrac{\partial G_2^0}{\partial \tilde{u}_{N/2}} & \dfrac{G_2^0}{\partial F} \\
        \vdots & \vdots & \ddots &  \vdots &  \vdots &   \vdots & \ddots & \vdots & \vdots \\
        \dfrac{\partial G_2^{N/2}}{\partial \tilde{\eta}_0} & \dfrac{\partial G_2^{N/2}}{\partial \tilde{\eta}_1} & \cdots & \dfrac{\partial G_2^{N/2}}{\partial \tilde{\eta}_{N/2}} & \dfrac{\partial G_2^{N/2}}{\partial \tilde{u}_0} & \dfrac{\partial G_2^{N/2}}{\partial \tilde{u}_1} & \cdots & \dfrac{\partial G_2^{N/2}}{\partial \tilde{u}_{N/2}} & \dfrac{G_2^{N/2}}{\partial F} \\
        \dfrac{\partial G_3}{\partial \tilde{\eta}_0} & \dfrac{\partial G_3}{\partial \tilde{\eta}_1} & \cdots & \dfrac{\partial G_3}{\partial \tilde{\eta}_{N/2}} & 0 & 0 & \cdots & 0 & 0
        \end{bmatrix}
    \end{align*}
    where the partial derivatives are approximated by a first-order finite difference scheme:
    \begin{align*}
        &\dfrac{\partial G}{\partial \tilde{\eta}_j} = \dfrac{G(\tilde{\eta}_0,\dots,\tilde{\eta}_j + \delta,\dots,\tilde{\eta}_{N/2},\tilde{u}_0,\dots,\tilde{u}_{N/2},F) - G(\tilde{\eta}_0,\dots,\tilde{\eta}_{N/2},\tilde{u}_0,\dots,\tilde{u}_{N/2},F)}{\delta}, \\
        &\dfrac{\partial G}{\partial \tilde{u}_j} = \dfrac{G(\tilde{\eta}_0,\dots,\tilde{\eta}_{N/2},\tilde{u}_0,\dots,\tilde{u}_j + \delta,\dots,\tilde{u}_{N/2},F) - G(\tilde{\eta}_0,\dots,\tilde{\eta}_{N/2},\tilde{u}_0,\dots,\tilde{u}_{N/2},F)}{\delta}, \\
        &\dfrac{\partial G}{\partial F} = \dfrac{G(\tilde{\eta}_0,\dots,\tilde{\eta}_{N/2},\tilde{u}_0,\dots,\tilde{u}_{N/2},F+\delta) - G(\tilde{\eta}_0,\dots,\tilde{\eta}_{N/2},\tilde{u}_0,\dots,\tilde{u}_{N/2},F)}{\delta},
    \end{align*}
    with $j = 0,\dots,N/2$ and $\delta = 10^{-10}$.
    
    The stopping criterion used is:
    \begin{equation*}
        \dfrac{\sum_{n = 1}^2 \sum_{j= 0}^{N/2} |G_n^j(\tilde{\eta}_0,\dots,\tilde{\eta}_{N/2},\tilde{u}_0,\dots,\tilde{u}_{N/2},F)| + |G_3(\tilde{\eta}_0,\dots,\tilde{\eta}_{N/2},\tilde{u}_0,\dots,\tilde{u}_{N/2},F)|}{N+3} < 10^{-10}.
    \end{equation*}
    
        To guarantee the convergence of Newton's method, the initial guess must be close to the solution. As a initial guess for $\tilde{\eta}$ we will use the analytical solution of the normalized KdV equation
    \begin{equation*}
    \eta(x) =  A\text{sech}^2(Kx), \quad K = \dfrac{1}{2}\sqrt{3A\dfrac{\alpha}{\beta}},\\
    \end{equation*}
    and as initial guess for $\tilde{u}$ we consider
    \begin{equation*}
        u(x) = \eta(x) - \dfrac{\alpha}{4}\eta^2(x) + \dfrac{\beta}{6}\eta_{xx}(x).
    \end{equation*}
    
    Finally, we used the negative of the  wave speed from the normalized KdV equation as an initial guess for $F$, i.e.,
    \begin{equation*}
        F = - \left( 1 + \dfrac{\alpha}{2} A \right).
    \end{equation*}  The numerical code  for the procedure presented here is  available at the appendix      \ref{AppendixA}.

    \section{Results and Discussions}
    
    In this section, we will compare the solution of the Boussinesq model with and without the sponge in three different contexts to validate and evaluate in which situations the model with the sponge is more advantageous. It is worth noting that close to the boundaries the models behave differently, since in the model with a sponge, the waves do not reach the boundaries, while in the model without a sponge, they do. Therefore, we evaluate the behavior of the models between the interval $[-80,80]$, which is the region where is expected that the sponge layer does not affect the solution. 
    
    Moreover, to verify the efficiency of the sponge layer, we study the relative error, given by 
    \begin{equation}
        E(t) = \dfrac{\Vert \vec{\eta}_c(t) - \vec{\eta}_s(t) \Vert_2}{\Vert \vec{\eta}_s(t)\Vert_2}, 
        \label{errorelativoparacompararmodeloscsesponja}
    \end{equation}
    where $\vec{\eta}_s$ and $\vec{\eta}_c$ are the vectors in which the coordinates are $\eta_{s_j} = \eta_s(x_j,t)$ and $\eta_{c_j} = \eta_c(x_j,t)$, where $\eta_s$ and $\eta_c$ are the solutions of the Boussinesq equations \eqref{boussinesqcomesponjaa2} with and without the sponge, respectively, and $x_j$ are common points in both grids such that $x_j \in (-80,80)$. Those grids are defined in the following way
    \begin{equation*}
        x_j = -L_s + j\Delta x, j = 0,1,\dots,N_s-1, L_s = \dfrac{\Delta x}{2}N_s,
    \end{equation*}
    and
    \begin{equation*}
        x_j = -L_c + j\Delta x, j = 0,1,\dots,N_c-1, L_s = \dfrac{\Delta x}{2}N_c.
    \end{equation*}

    \subsection{Solitary traveling waves}
    
\begin{figure}[!htb]
        \centering
    
        \includegraphics{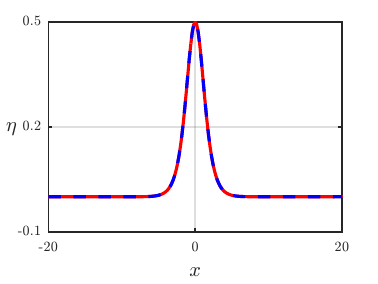}
        \includegraphics{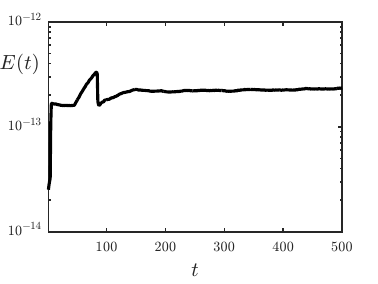}
        \caption{Left: Comparison between the surface profile of the Boussinesq model without sponge (in red) and with sponge (in blue). Right: Relative error graph.}
        \label{boussinesq-hp0}
    \end{figure}

    We consider a scenario where we have a traveling solitary wave in a channel with an even bottom and a constant pressure on the free surface. The initial conditions were the stationary solution of the Boussinesq equations, and for the numerical simulations, we fixed $N_c = N_s = 2^{10}$, $\Delta x = 0.2$, $\alpha = \beta = 0.01$ and $\Delta t = 0.01$. Figure~\ref{boussinesq-hp0} (left graph) displays the traveling waves at time $t = 352$, which is the time needed for the wave to travel 10 times its effective length. We choose the Froude number so that the wave is stationary, that is, $F = -c$, where $c$ is the wave speed. 
    It can be seen that the solutions practically coincide, which is reinforced by the relative error, which is of the order of $10^{-12}$, as shown in figure \ref{boussinesq-hp0} (right plot).  
    
    The low order of the relative error indicates that the sponge layer does not affect the region of the solution where $s(x) \approx 0$. It also shows that the choice of the sponge is appropriate and that the formulation is not numerically unstable.


    

    \subsection{Trapped waves}
    
    In this section, we study trapped waves in the Boussinesq equations framework within the sponge layer method. Trapped waves are waves that remain trapped in a certain region of space. They are characterized as waves that are temporarily or permanently trapped in a region of space, such as regions of low pressure or bounded by two topographical obstacles \cite{Flamarion:2021}. Here, we study the second case, considering a topography with two obstacles described by 
    \begin{equation}
     h_0(x) = A_0\left(e^{-(x+B_0)^2} + e^{-(x-B_0)^2}\right),
     \label{topografiaparaondaspresas}
    \end{equation}
    where the parameter $A_0$ controls the amplitude of the obstacles and $B_0$ denotes the point at which they are located. 
    
\begin{figure}[!h]
        \centering 
       \includegraphics{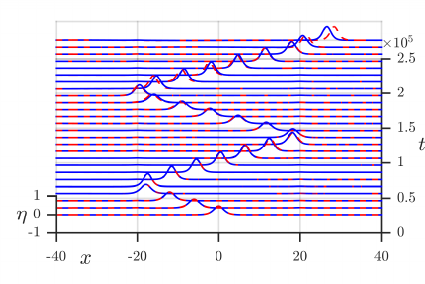} 
        \includegraphics{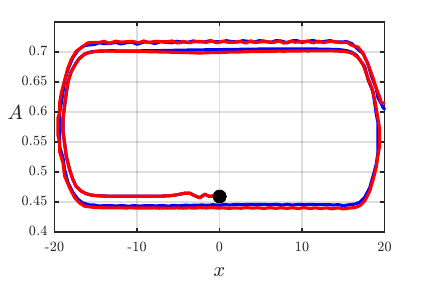}
        
        \caption{Left: Trapped waves between two bumps generated by Boussinesq equations with sponge (in blue)
    and without sponge (in red). Right: Amplitude of the waves generated by Boussinesq equations with sponge
    (in blue) and without sponge (in red).}
        \label{Boussinesq-mesh-op}
    \end{figure}

    Figure \ref{Boussinesq-mesh-op} (left) shows the evolution of a trapped wave for the Boussinesq model with and without sponge. It can be seen that both waves remain close for a period of time, but as time elapses, the wave solution of the Boussinesq model without sponge starts moving slightly faster than the wave of the model with sponge. This becomes more evident after they escape the region limited by the two obstacles. Moreover, we can observe (figure \ref{Boussinesq-mesh-op} right) that when the wave hits
    the obstacle located upstream, it increases its amplitude, while the opposite happens when the wave interacts
    with the obstacle located downstream.  The parameters used in this simulation are:  amplitude of the initial wave $A = 0.46$, $F = F_0 - 0.06\varepsilon$, where $F_0$ is the negative speed of the solitary traveling wave for a flat bottom computed trough the numerical method presented in the previous section,   $N_c = 2^{10}, N_s = 2^{13}, \Delta x = 0.2, \varepsilon = \alpha = \beta = 0.01, A_0 = \varepsilon^2, B_0 = 20$, and $\Delta t = 0.01$. 
    
\begin{figure}[!h]
        \centering
        
        \includegraphics{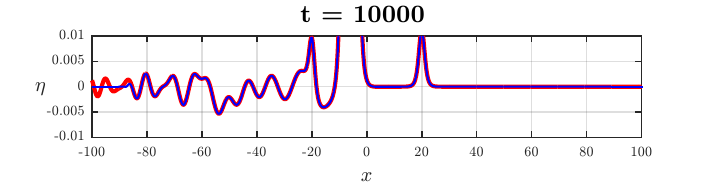}
        \includegraphics{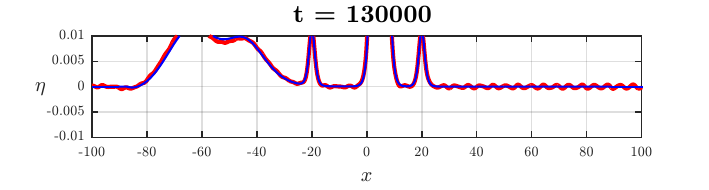}
        \includegraphics{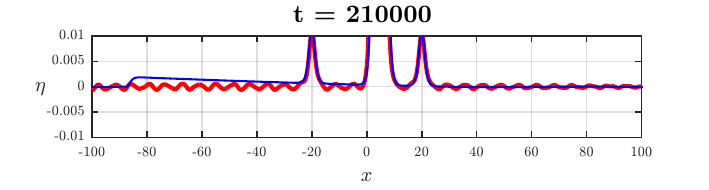}
        \caption{"Zoom in" of the trapped waves generated by Boussinesq equations with sponge (in blue) and without sponge (in red) for three time instants.}
        \label{snapshots-ondaspresas-boussinesq}
    \end{figure}

\begin{figure}[!h]
        \centering
        \includegraphics{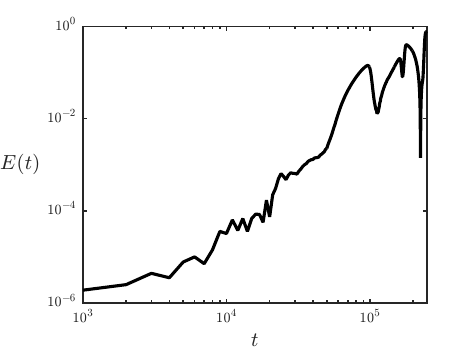}
        \caption{Relative error graph.}
        \label{Boussinesq-op}
    \end{figure}

    It is worth mentioning that there is a limitation in the choice of the length of the spatial domain in the model without sponge. As mentioned before, the numerical method used to solve the equations imposes periodic boundary conditions, and as consequence, when the solution reaches one of the boundaries, it re-enters from the other and start interacting with the solution. Because of this, and the fact that it takes a long period of time until the wave escapes the region, it would be necessary to choose a larger spatial domain. However, due to the computational costs, we could not choose more than $2^{13}$ Fourier modes. Figure \ref{snapshots-ondaspresas-boussinesq} shows the ``zoom in'' of the surface wave graph of the Boussinesq model with and without sponge for three different times. Note that the interaction between the trapped wave and the topographical obstacles generate waves of small amplitude $(t = 10000)$, which travel upstream. Since the domain is periodic, in the model without sponge, these waves go round the domain and start interacting with the trapped wave $(t = 130000)$. This spurious interaction ends up affecting the model without sponge $(t = 210000)$, which begins to differ from what is observed in the model with sponge.  The discrepancy between the waves is reflected in the error $E(t)$, which grows over time and becomes close to $10^0$, as shown in figure \ref{Boussinesq-op}. We observe the error increases rapidly when the time reaches order $10^{4}$, which indicates that the source of the error may come from the spurious interaction and not from the model.
    
    Hence, we can conclude that the Boussinesq model with sponge is advantageous for studying trapped waves, as the sponge is able to absorb small amplitude waves. In addition, this model manages to capture the same dynamics of trapped waves as the model without sponge, with a smaller spatial grid.
    
    
    \subsection{Waves generated by a current-topography interaction}
    
    We now examine a scenario where we have waves generated by an interaction of a current with a topography, under the effect of a constant pressure on the free surface. In this case we consider the initial conditions to be zero, and for the  first experiment we set $N_s = N_c = 2^{12}, \Delta x = 0.1, \varepsilon = \alpha = \beta = 0.01$, and $\Delta t = 0.01$. Since we are interested in evaluating the performance of the model with sponge rather than studying the surface profile generated, we fixed the topography as $h(x) = \varepsilon \frac{e^{-x^2}}{\sqrt{\pi}}$, the pressure as $P(x,t) = 0$ and $F = 1$. 
    
\begin{figure}[!htb]
    \centering
    \includegraphics{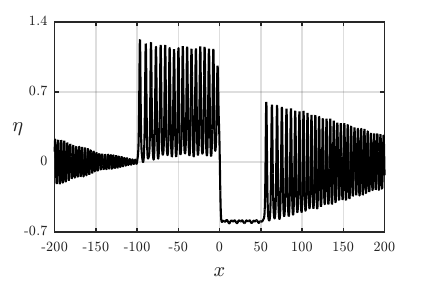}
    \includegraphics{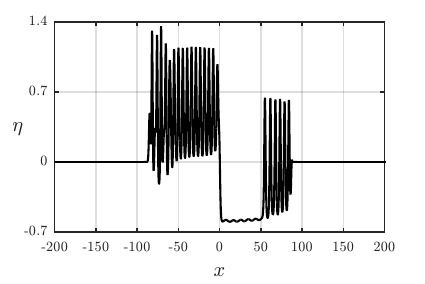}
    \caption{Evolution of the free surface at time $t = 16800$. Left: numerical solution without a sponge. Right: numerical solution with a sponge.}
    \label{boussinesq-tempo-demais}
    \end{figure}

        Figure \ref{boussinesq-tempo-demais} illustrates the solutions of the Boussinesq model with and without a sponge at time $t = 16800$. We observe that in the absence of a sponge (left plot), the downstream waves reach the right boundary, re-enter from the left boundary and start interacting with the upstream waves. Although this does not happen in the presence of a sponge (right plot), it is not capable of absorbing the high amplitude upstream waves. 
    
\begin{figure}[!htb]
        \centering
        \includegraphics{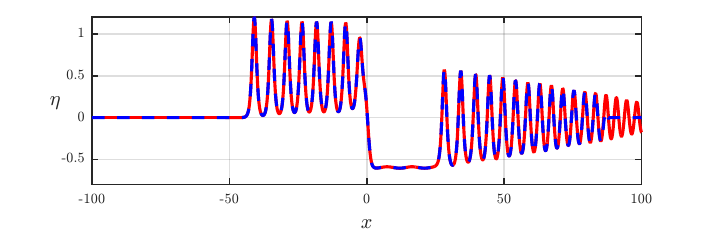}
        \includegraphics{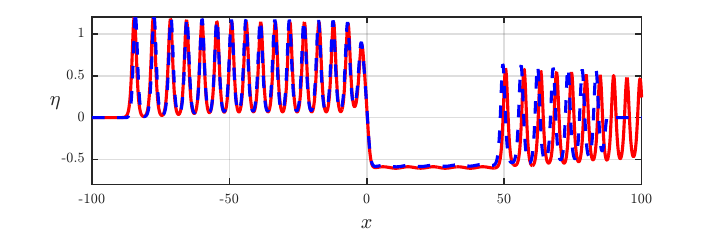}
        \caption{Comparison between the surface profile generated by Boussinesq equations with sponge (in blue) and without sponge (in red) at time $t = 7500$ (top) and $t = 14800$ (bottom).}
    \label{Boussinesq-og-215e212}
    \end{figure}

\begin{figure}[!htb]
        \centering
        \includegraphics{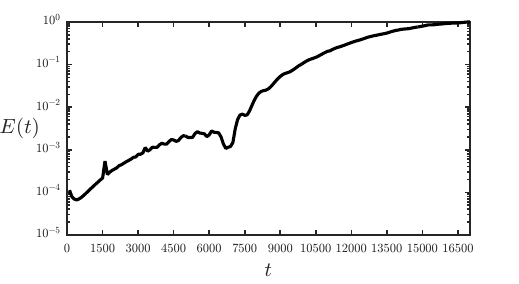}
        \caption{Relative error graph $E(t)$.}
        \label{errorelativo8}
    \end{figure}

    Given that we want to compare the models for as long as possible, we had to choose a larger spatial domain for the Boussinesq model without sponge, therefore, we increase the computational domain length by increasing the number of points for the numerical simulations for $N_s = 2^{15}$. As we can observe in figure \ref{Boussinesq-og-215e212}, the waves from both Boussinesq with and without a sponge remain close for a certain period of time. After that, those generated by the model without a sponge are slightly faster, particularly those going downstream. This difference impacts the relative error, which gets close of the order $10^0$, as shown in figure \ref{errorelativo8}
    .
    Furthermore, we notice the error starts increasing at time $t = 7500$, when the waves begin to display differences. After $t = 150000$, the error is mostly affected by the model, since it is not able to absorb waves with high amplitudes.
    
\begin{figure}[!htb]
        \centering
        \includegraphics{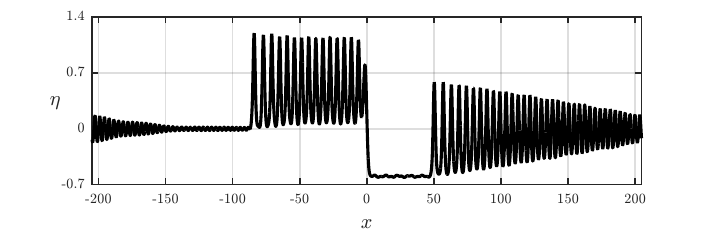}
        \caption{Surface profile generated by Boussinesq model without sponge with $N_s = 2^{12}$ at time $t = 14800$.}
        \label{boussinesq-tempo-malhamenor}
    \end{figure}

    Although the relative error graph indicates that there is a discrepancy between the models, the snapshot in figure \ref{Boussinesq-og-215e212} (bottom plot) shows that is not because of a difference between the solutions captured by the models, instead, it is due to the difference of positions of the waves $\eta_s$ and $\eta_c$. Because of this, it is possible to conclude that the Boussinesq model with a sponge captures the same dynamic as the model without, requiring a smaller spatial domain. In addition, it allows to simulate for long periods of time, because if we try to simulate without sponge at time $t = 14800$, with the same size of spatial grid as the model with sponge, we  see that the waves going downstream reach the boundaries, re-enter and interact with those going upstream, as showed in figure \ref{boussinesq-tempo-malhamenor}. However, instead of absorbing high amplitude waves going upstream, the model deflect them.  
    
    \section{Conclusion}
    
    The aim of this study was to propose a Boussinesq model with a sponge layer, and to study its advantages and limitations. Firstly, we validated the model by computing solitary traveling waves of the model with and without the sponge, in a scenario where both topography and the pressure on the free surface were constant. We found that there are no differences between these models in the region of the domain where $s(x) \approx 0$, leading us to conclude that the function adopted function to model the sponge layer  is adequate, as well as assuring us that the formulation with the sponge is not numerically unstable. 
    
    Afterwards, we compared the solutions of those models in the following contexts: (i) trapped waves; and (ii) waves generated by a current-topography interaction. Overall, we concluded that the model with sponge is able to capture the same dynamics of the waves generated by the model without the sponge, with a smaller computational domain, and allows numerical experiments to be carried out for a longer time. In summary, the model with sponge is useful because it allows to reduce the computational cost for the numerical simulations.

    \section*{Acknowledgments}
    The work of  L.G.M  was financed in part by the Coordena\c c\~ao de Aperfei\c coamento de Pessoal de N\'ivel Superior  Brasil (CAPES) -- Finance Code 001. The work of M.V.F. and R.R.Jr was supported in part by   National Council Scientific and Technological Development (CNPq) under  Chamada  CNPq/MCTI/N$^\circ$  10/2023-Universal.

    \section*{Data Availability Statement}
    
    No data was used for the research described in the article.
    
    \appendix
    
    \section{Numerical code}

    In this appendix  we present the MATLAB codes used   in the paper.

    \subsection{The stationary Boussinesq equations }  \label{AppendixA}

    Matlab code for  equation  \eqref{esquema}.

    \begin{verbatim}
    %  Parameters defined by the user:
    %   a:  parameter of nonlinearity (alpha)
    %   b:  dispersion parameter (beta) 
    %   A: amplitude of the wave 
    
    %   This code returns:
    %  eta: traveling wave profile
    %  u: horizontal velocity
    %  F
    
    % constant
    A = 0.44; 
    epsilon = 0.01;
    a = epsilon; 
    b = epsilon; 
    
    
    % Spatial domain
    n = 2^10;
    dx = 0.2;
    x = (-n/2:1:n/2-1)*dx;
    
    % Fourier modes
    dk = (2*pi)/(n*dx);
    k = [0:1:n/2-1 0 -n/2+1:1:-1]*dk;
    Mk = 1i*k;
    
    %%%%%%%.    Newton's method.    %%%%%%%%%
    K = 0.5*sqrt(3*(a/b)*A); 
    ddeta = -2*A*K^2*sech(K*x).^4 + 4*A*K^2*sech(K*x).^2.*tanh(K*x).^2;
    
    % Initial guess
    eta = A*sech(K*x).^2; 
    u = eta-0.25*a*eta.^2 + (b/6)*ddeta;
    F = -(0.5*a*A+1);
    
    ETA = fft(eta);
    U = fft(u);
    
    % Newton's method parameters
    delta = 10^-10; 
    Nimax = 100;  
    tol = 10^-10;
    
    Ni = 0;
    JacG = zeros(n+3); 
    
    Error = 1;
    while Error > tol && Ni <= Nimax
        g = G(ETA,U,F,Mk,A,a,b,n);
        
        % Compute the derivatives with respect to ETA
        ETAdelta = ETA;
        ETAdelta(1) = ETA(1) + delta;
        Gdelta = G(ETAdelta,U,F,Mk,A,a,b,n);
        JacG(:,1) = (Gdelta-g)./delta;
        
        for j = 2:1:n 
            ETAdelta = ETA;
            ETAdelta(j) = ETA(j) + delta;
            ETAdelta(n-j+2) = ETA(n-j+2) + delta;
            Gdelta = G(ETAdelta,U,F,Mk,A,a,b,n);
            JacG(:,j) = (Gdelta-g)./delta;
        end
    
        ETAdelta = ETA;
        ETAdelta(n/2+1) = ETA(n/2+1) + delta;
        Gdelta = G(ETAdelta,U,F,Mk,A,a,b,n);
        JacG(:,n/2+1) = (Gdelta-g)./delta;
        
        % Compute the derivatives with respect to U
        Udelta = U;
        Udelta(1) = U(1) + delta;
        Gdelta = G(ETA,Udelta,F,Mk,A,a,b,n);
        JacG(:,n/2+2) = (Gdelta-g)./delta;
        
         for j = 2:1:n/2
            Udelta = U;
            Udelta(j) = U(j) + delta;
            Udelta(n-j+2) = U(n-j+2) + delta;
            Gdelta = G(ETA,Udelta,F,Mk,A,a,b,n);
            JacG(:,n/2+j+1) = (Gdelta-g)./delta;
        end
        
        Udelta = U;
        Udelta(n/2+1) = U(n/2+1) + delta;
        Gdelta = G(ETA,Udelta,F,Mk,A,a,b,n);
        JacG(:,n+2) = (Gdelta-g)./delta;
        
        % Compute the derivative with respect to F
         Fdelta = F + delta;
         Gdelta = G(ETA,U,Fdelta,Mk,A,a,b,n);
         JacG(:,n+3) = (Gdelta-g)./delta;
        
          v=-JacG\g;    
         
         ETA =ETA+[v(1:n/2+1); v(n/2:-1:2)].'; 
         U=U+[v(n/2+2:n+2); v(n+1:-1:n/2+3)].';
         F = F + real(v(n+3));
    
         Error = sum(abs(g))/(n+3);
         
         clc
    
         Ni = Ni + 1
         display(['Error = ', num2str(Error)])
    end
    
      if Ni == Nimax
        display('The number of iterations has been exceeded')
       end
    
    % Solution
    eta = real(ifft(ETA));
    u = real(ifft(U));
    F;
    
    save('Stationary_solution.mat',...
        'epsilon','a','b','F','n','dx','eta','u')
    
    
    \end{verbatim}
    
    The function $G$, which is used in the loop, code is defined by the following way: 
    \begin{verbatim}
    function z = G(ETA,U,F,Mk,A,a,b,n)
    eta = real(ifft(ETA));
    u = real(ifft(U));
    
    
    G1 = F*ETA + U + a*fft(eta.*u);
    G2 = ETA + F*U + 0.5*a*fft((u).^2) - b*(F/3)*(Mk.^2).*U;
    G3 = eta(n/2+1)-A;
    
    z = [G1(1:1:n/2+1),G2(1:1:n/2+1),G3].';
    end 
    \end{verbatim}

    \subsection{The time-dependent Boussinesq equations} \label{AppendixB}

    Matlab code for  equation  \eqref{eqboussinesqnova5}.
    
    \begin{verbatim}
    %  Parameters and datas defined by the user:
    %   eta0:  initial wave profile
    %   u0:   initial horizontal velocity
    %   F : Froude number
    
    %   This code returns:
    %  eta: wave profile  for different times 
    %  u: horizontal velocity for different times
    
    
    % Here we are using the stationary solution computed in the previous code as  initial data
    load('Stationary_solution.mat',...
        'epsilon','a','b','F','n','dx','eta','u')
    
    eta0 = eta;
    u0= u;
    F0 = F;
    Nx = n;
    clear 'eta' 'u' 'F' 'n'
    
    x = (-Nx/2:1:Nx/2-1)*dx;
    
    % Fourier modes
    dk = (2*pi)/(Nx*dx);
    k = [0:1:Nx/2-1 0 -Nx/2+1:1:-1]*dk;
    Mk = 1i*k;
    
    % Topography
    C0 = 0.04;
    F = F0 - epsilon*C0;
    B = 20;
    h = 1*(epsilon^2)*(exp(-1*(x-B).^2) + exp(-1*(x+B).^2));
    H = fft(h);
    
    % Pressure
    p = 0*exp(-x.^2)./sqrt(pi);
    P = fft(p);
    
    % Time evolution
    dt = 0.01;
    taux = 3000/epsilon;
    Tf = taux; Ts = 10/epsilon; Ns = Tf/Ts+1; Nt = Ts/dt; 
    t = 0:Ts:Tf;       
    
    % Sponge
    aes = 10;
    x1 = -90; x2 = 90;
    s = 1*(0.5*aes*(tanh(x-x1)-tanh(x-x2))-aes);
    
    % Solution
    eta = zeros(Ns,Nx);  
    u= zeros(Ns,Nx);
    
    eta(1,:) = eta0;
    u(1,:) = u0;
    
    ETA = fft(eta(1,:));
    U =  fft(u(1,:));
    
    % Runge-Kutta (RK4)
    m = 1
    for ts = Ts:Ts:Tf
        
        for j = 1:Nt
            
            [K1,B1] =  camp(ETA,U,F,H,P,Mk,a,b,s,Nx);
            [K2,B2] =  camp(ETA+0.5*dt*K1,U+0.5*dt*B1,F,H,P,Mk,a,b,s,Nx);
            [K3,B3] =  camp(ETA+0.5*dt*K2,U+0.5*dt*B2,F,H,P,Mk,a,b,s,Nx);
            [K4,B4] =  camp(ETA+dt*K3,U+dt*B3,F,H,P,Mk,a,b,s,Nx);
            
            ETA = ETA + (dt/6)*(K1 + 2*K2 + 2*K3 + K4);
            U = U + (dt/6)*(B1 + 2*B2 + 2*B3 + B4);
            
        end
        eta(m+1,:) = real(ifft(ETA));
        u(m+1,:) = real(ifft(U));
        
        clc
        m = m+1 
    end
    \end{verbatim}
    
    The function $\text{camp}$, which is used in the loop, code is defined in the following way:
    \begin{verbatim}
    function [K1,B1] = camp(ETA,U,F,H,P,Mk,a,b,s,Nx)
    
    eta = real(ifft(ETA));
    u = real(ifft(U));
    h = real(ifft(H));
    
    ghat = fft((s.^2).*real(ifft(ETA)))./Mk;
    ghat(1) = 0; ghat(Nx/2+1) = 0;
    
    KK = -F*Mk.*ETA -Mk.*U +F*Mk.*H -a*Mk.*fft(eta.*u) +a.*Mk.*fft(h.*u)...
           + fft(2*s.*real(ifft(ETA)));
    BB = (-F*Mk.*U -Mk.*ETA -0.5*a*Mk.*fft(u.^2)+...
          (b/3)*F.*(Mk.^3).*U...
            -(b/2)*(F^2)*(Mk.^3).*H  - a*Mk.*P + ghat)./(1-(b/3).*(Mk.^2));
    
    K1 = KK;
    B1 = BB;
    end
    \end{verbatim}


\begin{thebibliography}{99}	
	
	

	
	 \bibitem{Alias:2019}
            Alias A, Ismail NNAN, Harun FN. Pseudospecteral method with linear damping effect and de-aliasing technique in solving nonlinear PDEs. {\it Journal of Physics: Conference Series}. 2019; 1366;012009.  http://doi.org/10.1088/1742-6596/1366/1/012009.
            
            	\bibitem{Baines}{ Baines P.}
    {\it Topographic effects in stratified flows.} { Cambridge: Cambridge University Press;} 1995.
    
    		
    	   
            \bibitem{Chen:2003}
             Chen M. Equations for bi-directional waves over an uneven bottom. {\it Mathematics and Computers in Simulation}. 2003;62:3-9. https://doi.org/10.1016/S0378-4754(02)00193-3.       
            
      
                      
                     
             
     
    
       
        \bibitem{Flamarion:2021}
        Flamarion MV, Ribeiro-Jr R. Trapped solitary-wave interaction for Euler equations with low-pressure region. {\it Computational and Applied Mathematics}. 2021;40:1-11. http://doi.org/10.1007/s40314-020-01407-0.
        
        
           \bibitem{Flamarion:sd}
            Flamarion MV, Ribeiro-Jr R, Vianna DLSS, et al. Trapped Solitary Waves in a Periodic External Force: A Numerical Investigation Using the Whitham Equation and the Sponge Layer Method. {\it Fluids}. 2023;8(8):223. https://doi.org/10.3390/fluids8080223.
    
           

        
        
           \bibitem{GrimShawMalewoong:2019}
            Grimshaw RHJ, Maleewong M. Transcritical flow over obstacles and holes: forced Korteweg-de Vries framework. {\it Journal of Fluid Mechanics}. 2019;881:660-678. http://doi.org/10.1017/jfm.2019.767       
            
            
             \bibitem{Grimshaw94}{Grimshaw R, Pelinovsky E, Tian X.} 
          {Interaction of a solitary wave with an external force.}
    	{\it Physica D.} 1994; 77: 405-433. http://doi.org/10.1002/sapm1996973235.
    	
    		\bibitem{Johnson} {Johnson RS.} 
          {Models for the formation of a critical layer in water wave propagation.}
    	{\it Phil. Trans. R. Soc. A.} 2012; 370:1638-1660. https://doi.org/10.1098/rsta.2011.0456.
    	
    	
    	\bibitem{Kim} {Kim H, Choi H.} 
          {A study of wave trapping between two obstacles in the forced Korteweg-de Vries equation.}
    	{\it J Eng Math.} 2018; 108:197-208. https://doi.org/10.1007/s10665-017-9919-5.
    	
    		\bibitem{Lee} {Lee S.} 
          {Dynamics of trapped solitary waves for the forced KdV equation.}
    	{\it Symmetry.} 2018; 10(5):129. https://doi.org/10.3390/sym10050129.
    	
    \bibitem{LeeWhang} {Lee S,  Whang S.} 
          { Trapped supercritical waves for the forced KdV equation with two bumps.}
    	{\it Appl Math Model.} 2015;  39:2649-2660. https://doi.org/10.1016/j.apm.2014.11.007.
	
	
	 \bibitem{Mani:2012}
            Mani A. Analysis and optimization of numerical sponge layers as a nonreflective boundary treatment. {\it Journal of Computational Physics}. 2012;231:704-716. http://doi.org/10.1016/j.jcp.2011.10.017

	
	
	  \bibitem{MejiaGrajales:2021}
            Mej'a LF, Grajales JCM. Analytical and Rothe time-discretization method for a Boussinesq-type system over an uneven bottom. {\it Communications in Nonlinear Science and Numerical Simulation}. 2021;102:105951. https://doi.org/10.1016/j.cnsns.2021.105951
            
              
              
           
    
    	\bibitem{Paul}{Milewski PA.} 
    	{The Forced Korteweg-de Vries equation as a model for waves generated by topography.}
    	{\it CUBO A mathematical Journal.} 2004; 6:33-51. 
    
            \bibitem{Pratt}{Pratt LJ. } 
    {On nonlinear flow with multiple obstructions.}
    {\it  J. Atmos. Sci.} 1984; 41:1214-1225. http://doi.org/10.1175/1520-0469(1984)041\%3C1214:ONFWMO\%3E2.0.CO;2.
    
    
    	  \bibitem{Shentangwang:2011}
            Shen J, Tang T, Wang L. {\it Spectral Methods: Algorithms, Analysis and Applications}. Berlin, BE: Springer; 2011.

   \bibitem{Trefethen:2001} 
        Trefethen LN. {\it Spectral Methods in Matlab}. Philadelphia: SIAM; 2001.
    	
    	
    	      \bibitem{Whitham:1974}
            Whitham GB. {\it Linear and Nonlinear Waves}. New York, NY: John Wiley \& Sons; 1974.
    
    
    \bibitem{Wu1}{Wu TY.} 
            {Generation of upstream advancing solitons by moving disturbances.}
    	{\it J Fluid Mech.} 1987;  184: 75-99. http://doi.org/10.1017/S0022112087002817.
    		
    		\bibitem{Wu2}{Wu DM, Wu TY.} 
           {Three-dimensional nonlinear long waves due to moving surface pressure. In: Proc. 14th. Symp.
    on Naval Hydrodynamics.}
    {\it Nat. Acad. Sci., Washington, DC.} 1982; 103-25.	
    		
    		
    		
    		
    	\end{thebibliography}
    	\end{document}